\documentclass[12pt]{article}
\usepackage[noend]{algpseudocode}
\usepackage{algorithmicx,algorithm}
\usepackage{amsmath}
\usepackage{indentfirst}
\usepackage{exscale}
\usepackage{relsize}
\usepackage{fancyhdr}  
\usepackage{amssymb}
\usepackage{mathrsfs}
\usepackage{color}
\usepackage{amsfonts}
 \usepackage{eucal}
\usepackage{bbding}
\allowdisplaybreaks[3]
\newtheorem{theorem}{Theorem}[section]
\newtheorem{lemma}[theorem]{Lemma}

\newtheorem{proposition}[theorem]{Proposition}
\newenvironment{proof}{{\noindent\it Proof.}}{\hfill $\square$\par}

\numberwithin{equation}{section}
  \usepackage{a4wide}

\def \O{{\Omega}}

\def\O{{\mathcal O}}

\begin{document}
\baselineskip 14pt
\title{\bf On an additive problem involving fractional powers with one prime and an almost prime variables}
\author{{Liuying  Wu}\\{ School of Mathematical Science, Tongji University}
\\ {Shanghai, 200092, P. R. China. }\\{Email: liuyingw@tongji.edu.cn}}

\date{}
\maketitle
\noindent {\bf Abstract}: For any real number $t$, let $[t]$ denote the integer part of $t$. In this paper it is proved that if $1<c<\frac{247}{238}$, then for sufficiently large integer $N$, the equation
\[\left[p^{c}\right]+\left[m^{c}\right]=N\]
has a solution in a prime $p$ and an almost prime $m$ with at most $\left[\frac{450}{247-238c}\right]+1$ prime factors. This result constitutes an improvement upon that of Petrov and Tolev \cite{PeTo}.



\section{\bf Introduction} \label{s1}

For fixed integer $k \geq 1$ and sufficiently large integer $N$, the well-known Waring problem is devoted to investigating the solvability of the following Diophantine equality
\begin{equation}\label{e.101}
N=m_{1}^{k}+m_{2}^{k}+\cdots+m_{s}^{k}
\end{equation}
in integer variables $m_{1}, m_{2}, \ldots, m_{s}$. In 1933, Segal \cite{Seg1,Seg} studied the following anolog of the equation \eqref{e.101}. Suppose that $c>1$ and $c\notin \mathbb{N}$, there exists a positive integer $s=s(c)$ such that for every sufficiently large natural number $N$, the equation
\begin{equation}\label{e.102}
N=\left[m_{1}^{c}\right]+\left[m_{2}^{c}\right]+\cdots+\left[m_{s}^{c}\right]
\end{equation}
has a solution with $m_{1}, m_{2}, \ldots, m_{s}$ integers, where $[t]$ denotes the integer part of any $t\in \mathbb{R}$. These sequences of the form
$$([n^c])_{n=1}^{\infty},\qquad c>1,\qquad c\notin \mathbb{N}$$
are so-called \emph{Piatetski-Shapiro sequence} in honour of Piatetski-Shapiro, who \cite{PS} showed that such sequence contains infinitely many prime numbers if $1<c<\frac{12}{11}.$ The range of $c$ has been improved many times and the best result till now is $1<c<\frac{243}{205}$ thanks to Rivat and Wu \cite{RW}.
\vskip 2 pt

For the special case $s=2$ in \eqref{e.102}, many mathematicians have derived many splendid results. In 1973, Deshouillers \cite{De} proved that if $1<c<\frac{4}{3}$, then for every sufficiently large integer $N$ the equation
\begin{equation}\label{e.103}
N=\left[m_{1}^{c}\right]+\left[m_{2}^{c}\right]
\end{equation}
has a solution with $m_{1}$ and $m_{2}$ integers. Later, the range of $c$ was enlarged to $\frac{55}{41}$ and $\frac{3}{2}$ by Gritsenko \cite{Gr} and Konyagin \cite{Kon}, respectively. On the other hand, Kumchev \cite{Kum} proved that if $1<c<\frac{16}{15}$, then for every sufficiently large integer $N$ can be represented in the form \eqref{e.103}, where $m_{1}$ is a prime and $m_{2}$ is an integer.
\vskip 2 pt
For any natural number $r$, let $\mathscr{P}_{r}$ denote an almost-prime with at most $r$ prime factors, counted according to multiplicity. The celebrated theorem of Chen \cite{Chen} states that every sufficiently large integer $N$ can be represented a sum of a prime and an almost prime $\mathscr{P}_{2}$. Bearing in mind this profound result, it is reasonable to conjecture that there exists a constant $c_{0}>1$ such that if $1<c<c_{0}$, then the equation \eqref{e.103} has a solution with $m_{1}$ a prime and $m_{2} \in \mathscr{P}_{2}$ for sufficiently large $N$. Motivated by Kumchev \cite{Kum}, Petrov and Tolev \cite{PeTo} proved that if $1<c<\frac{29}{28}$, then every sufficiently large integer $N$ can be represented as
\begin{equation}\label{e.104}
N=\left[p^{c}\right]+\left[m^{c}\right],
\end{equation}
where $p$ is a prime and $m$ is an almost prime with at most $\left[\frac{52}{29-28c}\right]+1$ prime factors.
\vskip 2 pt
In this paper, motivated by \cite{PeTo}, we shall prove the following sharper result:
\begin{theorem}
Suppose that $1<c<\frac{247}{238}$. Then every sufficiently large integer $N$ can be represented as
$$N=\left[p^{c}\right]+\left[m^{c}\right],$$
where $p$ is a prime and $m$ is an almost prime with at most $\left[\frac{450}{247-238c}\right]+1$ prime factors.
\end{theorem}
{\bf Remark} In order to compare our result with the result of \cite{PeTo}, we list the numerical result as follows:
\[\frac{29}{28}=1.03571 \cdots, \qquad \frac{247}{238}=1.03782 \cdots\]
Moreover, it is easy to verify that when $c\rightarrow 1$, the result of \cite{PeTo} indicates that $m \in \mathscr{P}_{53}$ in \eqref{e.104}, while our result shows that $m$ is $\mathscr{P}_{51}$. It is worth mentioning that by inserting a weighted sieve approach into our argument, one may further reduce the number of prime factors of $m$, but that is not the aim of this paper.

\section{\bf Preliminaries} \label{s2}
Throughout this paper, the letter $p$ and $q$ always stand for prime numbers. We use $\varepsilon$ to denote a sufficiently small positive number, and the value of $\varepsilon$ may change from statement to statement. For any natural number $n$, we use $\mu(n), \Lambda(n)$ and $\tau(n)$ to denote Möbius' function, von Mangolds' function and Dirichlet divisor function, respectively. We write $f=\mathcal{O}(g)$ or, equivalently, $f \ll g$ if $|f|\leq Cg$ for some positive number $C$. If we have simultaneously $A\ll B$ and $B\ll A$, then we shall write $A\asymp B$. Let $\{t\}$ be the fractional part of $t$, the function $\psi(t)$ is defined by $\psi(t)=\frac{1}{2}-\{t\}$. And we use $e(\alpha)$ to denote $e^{2\pi i\alpha}$. In addition, we define
\begin{equation}\label{e.2000}
1<c<\frac{247}{238}, \qquad \gamma=\frac{1}{c}, \qquad P=10^{-9} N^{\gamma}, \qquad \delta=\frac{247 \gamma-238}{225}-\varepsilon,
\end{equation}
\begin{equation}\label{e.200}
D=N^{\delta}, \qquad z=N^{\frac{\delta}{2}-\varepsilon}, \qquad P(z)=\prod_{2<p<z} p.
\end{equation}

\begin{lemma}
Let $r\geq1$ be an integer, and let $\alpha,\beta$ and $\Delta$ be real numbers such that $0<\Delta<\frac{1}{4}$ and $\Delta\leq \beta-\alpha\leq1-\Delta.$ Then there exists a function $\theta(x)$ that is periodic with period 1 and satisfies the conditions
\begin{itemize}
  \item[(1)] $\theta(x)=1$ on the interval $\alpha+\frac{\Delta}{2}\leq x\leq \beta-\frac{\Delta}{2};$
  \item[(2)] $0<\theta(x)<1$ on the intervals
  $$\alpha-\frac{\Delta}{2}<x<\alpha+\frac{\Delta}{2},\quad \text{and}\quad \beta-\frac{\Delta}{2}<x<\beta+\frac{\Delta}{2};$$
  \item[(3)] $\theta(x)=0$ on the interval $\beta+\frac{\Delta}{2}\leq x\leq1+\alpha-\frac{\Delta}{2};$
  \item[(4)] The Fourier expansion of $\theta(x)$ is of the form
  $$\theta(x)=\beta-\alpha+\sum_{\substack{m=-\infty\\ m\neq0}}^{+\infty}g(m)e(mx),$$
where
$$|g(m)|\leq\min\left(\beta-\alpha,\frac{1}{\pi|m|},\frac{1}{\pi|m|}\left(\frac{r}{\pi|m|\Delta}\right)^r\right).$$
\end{itemize}
\end{lemma}
\begin{proof}
See \cite[Chapter 1, Lemma A]{Kara}.
\end{proof}

\begin{lemma}
For any $H\geq1,$ there exist numbers $a_{h}$ and $b_{h}$ such that
$$\left|\psi(x)-\sum_{0<|h|\leq H}a_{h}e(hx)\right|\leq\sum_{|h|\leq H}b_{h}e(hx),\qquad a_{h}\ll\frac{1}{|h|},\qquad b_{h}\ll \frac{1}{H}.$$
\end{lemma}
\begin{proof}
See Vaaler \cite{Vaa}.
\end{proof}

\begin{lemma}
Suppose that $D>4$ is a real number and let $\lambda^{\pm}(d)$ be the Rosser's functions of level $D$. Then we have the following properties:
\begin{itemize}
  \item[(1)] For any positive integer d we have
\[\left|\lambda^{\pm}(d)\right|\leq 1, \quad \lambda^{\pm}(d)=0 \quad \text {if} \quad d>D \quad \text {or} \quad \mu(d)=0 .\]
  \item[(2)] If $n\in\mathbb{N}$ then
\[\sum_{d\mid n}\lambda^{-}(d)\leq\sum_{d \mid n}\mu(d)\leq\sum_{d\mid n}\lambda^{+}(d)\].
  \item[(3)] If $z\in\mathbb{R}$ is such that $z^{2}\leq D\leq z^{3}$ and if
\[P(z)=\prod_{2<p<z} p, \quad \mathscr{B}=\prod_{2<p<z}\left(1-\frac{1}{p}\right), \quad \mathscr{N}^{\pm}=\sum_{d\mid P(z)} \frac{\lambda^{\pm}(d)}{d}, \quad s_{0}=\frac{\log D}{\log z},\]
then we have
\begin{equation*}\label{e.201}
\mathscr{B}\leq\mathscr{N}^{+}\leq\mathscr{B}\left(F\left(s_{0}\right)+\mathcal{O}\left((\log D)^{-1/3}\right)\right),
\end{equation*}
\begin{equation}\label{e.202}
\mathscr{B}\geq\mathscr{N}^{-}\geq\mathscr{B}\left(f\left(s_{0}\right)+\mathcal{O}\left((\log D)^{-1/3}\right)\right),
\end{equation}
where
\begin{equation}\label{e.02}
f(s)=\frac{2e^G}{s}\log(s-1),\quad F(s)=\frac{2e^G}{s}\quad \text{for}~~2\leq s\leq3,
\end{equation}
and $G$ stands for Euler constant.
\end{itemize}
\end{lemma}
\begin{proof}
This is a special case of a more general case, one can see Greaves \cite{Gre}.
\end{proof}

\begin{lemma}
Suppose that $f(x):[a, b]\rightarrow\mathbb{R}$ has continuous derivatives of arbitrary order on $[a,b]$, where $1\leq a<b\leq2a$. Suppose further that
\[\left|f^{j}(x)\right|\asymp\lambda_{1}a^{1-j}, \qquad j\geq 1, \qquad x \in[a, b].\]
Then for any exponent pair $(\kappa,\lambda)$, we have
\[\sum_{a<n\leq b}e(f(n))\ll\lambda_{1}^{\kappa}a^{\lambda}+\lambda_{1}^{-1}.\]
\end{lemma}
\begin{proof}
See (3.3.4) of Graham and Kolesnik \cite{GK}.
\end{proof}

\begin{lemma}
For any complex numbers $z_{n}$, we have
\[\left|\sum_{a<n\leq b}z_{n}\right|\leq\left(1+\frac{b-a}{Q}\right)\sum_{|q|<Q}\left(1-\frac{|q|}{Q}\right)\sum_{a<n, n+q\leq b}z_{n+q}\overline{z_{n}},\]
where $Q$ is any positive integer.
\end{lemma}
\begin{proof}
See \cite[Lemma 8.17]{IK}.
\end{proof}

\section{\bf Outline of the method}\label{s3}
Let $N$ be a sufficiently large integer. A key point in our paper is the study of the sum
\begin{equation}\label{e.301}
\Gamma=\sum_{\substack{P<p\leq 2P, m\in\mathbb{N}\\ \left[p^{c}\right]+\left[m^{c}\right]=N\\ (m,P(z))=1}}(\log p).
\end{equation}
\vskip 2 pt
Now we consider $\Gamma.$ If $\Gamma>0,$ then there is a prime $p$ and a natural number $m$ satisfying
\begin{equation}\label{e.302}
\left[p^{c}\right]+\left[m^{c}\right]=N, \qquad (m,P(z))=1.
\end{equation}
It follows from \eqref{e.302} that any prime factor of $m$ is greater or equal to $z$. Suppose that $m$ has $l$ prime factors, counted with the multiplicity. Then by \eqref{e.200} we have
$$N^{(\frac{\delta}{2}-\varepsilon)l}=z^l\leq m\leq N^\gamma,$$
and thus $l\leq \frac{2\gamma}{\delta-2\varepsilon}.$ This implies that if $\Gamma>0$ then \eqref{e.104} has a solution with $p$ a prime and $m$ an almost prime with at most $\left[\frac{2\gamma}{\delta-2\varepsilon}\right]$ prime factors, we choose $\varepsilon$ small enough such that
\begin{equation}\label{e.303}
2<\frac{\delta}{\frac{\delta}{2}-\varepsilon}<3.
\end{equation}
With the choice \eqref{e.2000} of $c$ and $\delta$, it is not difficult to see that
$$\left[\frac{2\gamma}{\delta-2\varepsilon}\right]\leq \left[\frac{450}{247-238c}\right]+1.$$
Thus, Theorem 1.1 will be proved if we can show that
\[\Gamma\gg\frac{N^{2 \gamma-1}}{\log N}.\]
\vskip 2 pt
To get the desired result, we need following proposition, which plays a central role in the proof of Theorem 1.1.
\begin{proposition}
Let
$$\frac{238}{247}<\gamma<1,\qquad \delta=\frac{247 \gamma-238}{225}-\varepsilon,$$
and $\lambda(d)$ be the Rosser's weights of level $D=N^\delta$. Define
\begin{equation}\label{e.308}
\Sigma_{j}=\sum_{d\leq D}\lambda(d)\sum_{P<p\leq 2P}(\log p) \psi\left(-\frac{1}{d}\left(N+j-\left[p^{c}\right]\right)^{\gamma}\right), \qquad j=0,1.
\end{equation}
Then we have
\begin{equation*}\label{e.307}
\Sigma_{j}\ll\frac{N^{2\gamma-1}}{(\log N)^{2}}, \qquad j=0,1.
\end{equation*}
\end{proposition}

\section{\bf Proof of Proposition 3.1}\label{s4}

\subsection{\bf The estimation of the sums $\Sigma_{0}$ and $\Sigma_{1}$}

Consider the sum $\Sigma_{j}$ defined in \eqref{e.308}. By Lemma 2.2, we can write
\begin{equation}\label{e.01}
\Sigma_{j}=\Sigma_{j}^{(1)}+\O\left(\Sigma_{j}^{(2)}\right),
\end{equation}
where
\begin{flalign*}
\Sigma_{j}^{(1)}&=\sum_{d\leq D}\lambda(d)\sum_{P<p\leq 2P}(\log p)\sum_{0<|h|\leq H}a_{h}e\left(-\frac{h}{d}\left(N+j-\left[p^{c}\right]\right)^{\gamma}\right),\\
\Sigma_{j}^{(2)}&=\sum_{d\leq D}\lambda(d)\sum_{P<p\leq 2P}(\log p)\sum_{|h|\leq H}b_{h}e\left(-\frac{h}{d}\left(N+j-\left[p^{c}\right]\right)^{\gamma}\right).
\end{flalign*}

Let
\begin{equation}\label{e.002}
W(v)=\sum_{P<p \leq 2 P}(\log p) e\left(v\left(N+j-\left[p^{c}\right]\right)^{\gamma}\right).
\end{equation}
We begin with the sum $\Sigma_{j}^{(1)}$. Changing the order of summation and using the upper bound $a_{h}\ll|h|^{-1}$ we get
\begin{equation}\label{e.03}
\Sigma_{j}^{(1)}=\sum_{d\leq D}\lambda(d)\sum_{0<|h|\leq H}a_{h}W\left(-\frac{h}{d}\right)\ll\sum_{d\leq D}\sum_{1\leq h\leq H}\frac{1}{h}\left|W\left(\frac{h}{d}\right)\right|.
\end{equation}
For the sum $\Sigma_{j}^{(2)}$, by \eqref{e.002} we have
\begin{align}
\Sigma_{j}^{(2)}&\ll\sum_{d\leq D}\sum_{P<p\leq 2P}(\log p)\sum_{|h|\leq H}b_{h} e\left(-\frac{h}{d}\left(N+j-\left[p^{c}\right]\right)^{\gamma}\right) \nonumber\\
&=\sum_{d\leq D}\sum_{|h|\leq H}b_{h} W\left(-\frac{h}{d}\right)\ll\sum_{d\leq D}\sum_{|h|\leq H} \frac{1}{H}\left|W\left(\frac{h}{d}\right)\right|\nonumber\\
&\ll \sum_{d\leq D}\frac{1}{H}|W(0)|+\sum_{d\leq D}\sum_{1\leq h\leq H} \frac{1}{H}\left|W\left(\frac{h}{d}\right)\right|.
\end{align}
By \eqref{e.2000} and Chebyshev's prime number theorem we find that $W(0)\asymp N^{\gamma}$. Let
\begin{equation}\label{e.04}
H=dN^{1-\gamma}(\log N)^{3}.
\end{equation}
Now, using \eqref{e.01}, \eqref{e.03}-\eqref{e.04} we obtain
\begin{equation}\label{e.05}
\Sigma_{j} \ll \frac{N^{2 \gamma-1}}{(\log N)^{2}}+\sum_{d \leq D} \sum_{h \leq H} \frac{1}{h}\left|W\left(\frac{h}{d}\right)\right|, \quad j=0,1.
\end{equation}

\subsection{\bf The evaluation of the sum $W(v)$}

Next, we shall consider the exponential sum $W(v)$ defined in \eqref{e.002}. Applying Lemma 2.1 with parameters
$$\alpha=-\frac{1}{4Z},\qquad \beta=\frac{1}{4Z},\qquad \Delta=\frac{1}{2Z},\qquad r=[\log N],$$
where 
\begin{equation}\label{e.11}
Z\asymp dN^{1-\gamma}(\log N)^7
\end{equation}
is an integer. Then $\theta(x)$ is periodic with period 1 satisfies
\begin{equation*}
\theta(0)=1;\quad 0<\theta(x)<1 \quad \text{for}\quad 0<|x|<\frac{1}{2Z};\quad \theta(x)=0\quad \text{for}\quad \frac{1}{2Z}\leq x\leq\frac{1}{2}.
\end{equation*}
Furthermore, the Fourier series of $\theta(x)$ is given by
$$\theta(x)=\frac{1}{2Z}+\sum_{\substack{m\in\mathbb{Z} \\ m\neq 0}}g(m)e(mx), \quad \text {with} \quad\left|g(m)\right|\leq \min\left(\frac{1}{2Z}, \frac{1}{|m|}\left(\frac{2Z[\log N]}{\pi|m|}\right)^{[\log N]}\right).$$
By \cite[(37) and (38)]{PeTo}, we have
\begin{equation*}\label{e.06}
\theta(x)=\sum_{|m|\leq Z(\log N)^{4}}g(m)e(mt)+O\left(N^{-\log\log N}\right),\quad \left|g(m)\right|\leq \frac{1}{2Z}.
\end{equation*}
Let
\[\theta_{z}(x)=\theta\left(x-\frac{z}{2Z}\right) \quad \text { for } \quad z=0,1,2, \ldots, 2Z-1.\]
Then from \cite[(43)]{PeTo} and \cite[(45)]{PeTo}, we find
\begin{flalign}
\theta_{z}(x)&=\sum_{|m|\leq Z(\log N)^{4}}g_{z}(m)e(mt)+O\left(N^{-\log\log N}\right),\quad \left|g_{z}(m)\right|\leq \frac{1}{2Z},\label{e.17}
\end{flalign}
and
\begin{equation}\label{e.007}
\sum_{z=0}^{2Z-1}\theta_{z}(x)=1\qquad \text{for~all}\quad x\in\mathbb{R}.
\end{equation}

\vskip 2 pt 
Now, it follows from \eqref{e.002} and \eqref{e.007} that
\begin{equation}\label{e.08}
W(v)=\sum_{P<p \leq 2P}(\log p) e\left(v\left(N+j-\left[p^{c}\right]\right)^{\gamma}\right)\sum_{z=0}^{2Z-1}\theta_{z}(p^c)=\sum_{z=0}^{2Z-1}W_{z}(v),
\end{equation}
where
\begin{equation}\label{e.09}
W_{z}(v)=\sum_{P<p\leq 2P}(\log p)\theta_{z}(p^c)e\left(v\left(N+j-\left[p^{c}\right]\right)^{\gamma}\right).
\end{equation}
By \cite[(48) and (49)]{PeTo}, \eqref{e.2000} and \eqref{e.11}, we can see that the contribution with $z=0$ to $W(v)$ is
\begin{equation}\label{e.10}
W_{0}(v)\ll (\log N)\left(\frac{N^\gamma}{Z}+N^{\frac{1}{2}}Z^{\frac{1}{2}}\log^6N\right)\ll \frac{N^{2\gamma-1}}{d\log^6N}.
\end{equation}
For the contribution of $W_{z}(v)$ with $z\neq 0,$ using \cite[(52)]{PeTo} and \eqref{e.11} we get
\begin{equation}\label{e.12}
W_{z}(v)=V_{z}(v)+\O\left(\frac{vN^{2\gamma-2}}{d\log^7N}\sum_{P<p\leq 2P}(\log p)\theta_{z}(p^c)\right),
\end{equation}
where
\begin{equation}\label{e.16}
V_{z}(v)=\sum_{P<p\leq 2P}(\log p)\theta_{z}(p^c)e\left(v\left(N+j-\left[p^{c}\right]+\frac{z}{2Z}\right)^{\gamma}\right).
\end{equation}
Thus, combining \eqref{e.08}, \eqref{e.10} and \eqref{e.12} we have
\begin{align}
W(v)&=\sum_{z=1}^{2Z-1}W_{z}(v)+W_{0}(v)\nonumber\\
&=\sum_{z=1}^{2Z-1}V_{z}(v)+\O\left(\frac{vN^{2\gamma-2}}{d\log^7N}\sum_{P<p\leq 2P}(\log p)\sum_{z=1}^{2Z-1}\theta_{z}(p^c)\right)+\O\left(\frac{N^{2\gamma-1}}{d\log^6N}\right)
\end{align}
Now, we use \eqref{e.2000}, \eqref{e.007}, \eqref{e.10} and Chebyshev's prime number theorem to derive that
\begin{align*}
\Xi&=\frac{vN^{2\gamma-2}}{d\log^7N}\sum_{P<p\leq 2P}(\log p)\sum_{z=1}^{2Z-1}\theta_{z}(p^c)\\
&=\frac{vN^{2\gamma-2}}{d\log^7N}\left(\sum_{P<p\leq 2P}(\log p)-\sum_{P<p\leq 2P}(\log p)\theta_{0}(p^c)\right)\\
&\ll \frac{vN^{3\gamma-2}}{d\log^7N},
\end{align*}
and therefore
\begin{equation}\label{e.13}
W(v)=\sum_{z=1}^{2Z-1}V_{z}(v)+\O\left(\frac{vN^{3\gamma-2}}{d\log^7N}\right)
+\O\left(\frac{N^{2\gamma-1}}{d\log^6N}\right).
\end{equation}
\vskip 2 pt
From now on, we assume that
\begin{equation}\label{e.14}
v=\frac{h}{d}, \quad \text{where}\quad 1\leq d\leq D,\quad 1\leq h\leq H.
\end{equation}
Then from \eqref{e.04} we can see that the error term in \eqref{e.13} is $ \O\left(N^{2\gamma-1}(d\log^4N)^{-1}\right),$ which means
\begin{equation}\label{e.15}
W(v)=\sum_{z=1}^{2Z-1}V_{z}(v)+\O\left(\frac{N^{2\gamma-1}}{d\log^4N}\right).
\end{equation}
\vskip 2 pt
Now, we consider the sum $V_{z}(v)$ defined in \eqref{e.16}, where $v$ satisfies \eqref{e.14}. By \eqref{e.11} and \eqref{e.17} we find that
\begin{align}
V_{z}(v)& =\sum_{P<p\leq 2P}(\log p)\left(\sum_{|r|\leq Z(\log N)^{4}}g_{z}(r)e(rp^c)\right) e\left(v\left(N+j-p^{c}+\frac{z}{2Z}\right)^{\gamma}\right)+O\left(N^{-10}\right) \nonumber\\
&=\sum_{|r|\leq Z(\log N)^{4}}g_{z}(r)U\left(N+j+\frac{z}{2Z},r,v\right)+O\left(N^{-10}\right)\nonumber\\
&\ll N^{-10}+\frac{1}{Z}\sum_{|r|\leq R}\sup_{T\in[N,N+2]}|U(T,r,v)|,\label{e.18}
\end{align}
where
\begin{equation}\label{e.19}
U(T,r,v)=\sum_{P<p\leq2P}(\log p)e\left(rp^{c}+v\left(T-p^{c}\right)^{\gamma}\right),
\end{equation}
\begin{equation}\label{e.20}
R=dN^{1-\gamma}(\log N)^{12}.
\end{equation}
Finally, combining \eqref{e.05}, \eqref{e.15} and \eqref{e.18}, we have
\begin{equation}\label{e.S01}
\left|\Sigma_{0}\right|+\left|\Sigma_{1}\right|\ll\frac{N^{2\gamma-1}}{(\log N)^{2}}+\sum_{d\leq D}\sum_{h\leq H} \frac{1}{h}\sum_{|r|\leq R}\sup_{T\in[N, N+2]}|U(T,r,v)|.
\end{equation}

\begin{lemma}
Let $f(n)$ be a complex valued function defined on $n \in(P,2P]$. Then we have
\[\sum_{P<n\leq 2P}\Lambda(n)f(n)=S_{1}-S_{2}-S_{3},\]
where
\begin{flalign*}
S_{1}&=\sum_{k\leq P^{\frac{1}{3}}}\mu(k)\sum_{\frac{P}{k}<\ell\leq\frac{2P}{k}}(\log \ell)f(k\ell), \\
S_{2}&=\sum_{k\leq P^{\frac{2}{3}}}c(k)\sum_{\frac{P}{k}<\ell\leq\frac{2P}{k}}f(k\ell), \\
S_{3}&=\sum_{P^{\frac{1}{3}}<k\leq P^{\frac{2}{3}}}a(k)\sum_{\frac{P}{k}<\ell\leq\frac{2P}{k}}\Lambda(\ell)f(k\ell), \\
\end{flalign*}
and where $a(k), c(k)$ are real numbers satisfying
\[|a(k)|\leq\tau(k), \qquad|c(k)|\leq\log k.\]
\end{lemma}
\begin{proof}
Can be found in \cite{RCV}.
\end{proof}
\vskip 2 pt

By \eqref{e.19} and Lemma 4.1 with $f(n)=e(rn^{c}+v\left(T-n^{c})^{\gamma}\right)$, we have
\begin{equation}\label{e.404}
U(T,r,v)=U_{1}-U_{2}-U_{3}+\mathcal{O}\left(N^{\frac{\gamma}{2}}\right),
\end{equation}
where
\begin{flalign}
U_{1}&=\sum_{k\leq P^{\frac{1}{3}}}\mu(k)\sum_{\frac{P}{k}<\ell\leq\frac{2P}{k}}(\log \ell)f(k\ell), \label{e.405}\\
U_{2}&=\sum_{k\leq P^{\frac{2}{3}}}c(k)\sum_{\frac{P}{k}<\ell\leq\frac{2P}{k}} f(k\ell), \label{e.406}\\
U_{3}&=\sum_{P^{\frac{1}{3}<k\leq P^{\frac{2}{3}}}}a(k)\sum_{\frac{P}{k}<\ell\leq\frac{2P}{k}}\Lambda(\ell)f(k\ell),\label{e.407}
\end{flalign}
and $|a(k)|\leq\tau(k), |c(k)|\leq\log k$. Then from \eqref{e.19}, \eqref{e.S01} and \eqref{e.404} we obtain
\begin{equation}\label{e.408}
\left|\Sigma_{0}\right|+\left|\Sigma_{1}\right|\ll\frac{N^{2\gamma-1}}{(\log N)^{2}}+\sum_{i=1}^{3}\Omega_{i},
\end{equation}
where
\begin{equation}\label{e.409}
\Omega_{i}=\sum_{d\leq D}\sum_{h\leq H}\frac{1}{h}\sum_{|r|\leq R}\sup_{T\in[N,N+2]}\left|U_{i}\right|.
\end{equation}
Moreover, we write $U_{2}=U_{2}^{(1)}+U_{2}^{(2)}$, where
\[U_{2}^{(1)}=\sum_{k\leq P^{\frac{1}{3}}}c(k)\sum_{\frac{P}{k}<\ell\leq\frac{2P}{k}}f(k\ell), \quad U_{2}^{(2)}=\sum_{P^{\frac{1}{3}<k\leq P^{\frac{2}{3}}}}c(k)\sum_{\frac{P}{k}<\ell\leq\frac{2P}{k}}f(k\ell),\]
then we have
\begin{equation}\label{e.410}
\Omega_{2}\leq\Omega_{2}^{(1)}+\Omega_{2}^{(2)},
\end{equation}
where
\[\Omega_{2}^{(1)}=\sum_{d\leq D}\sum_{h\leq H}\frac{1}{h}\sum_{|r|\leq R}\sup_{T\in[N,N+2]}\left|U_{2}^{(1)}\right|, \quad \Omega_{2}^{(2)}=\sum_{d\leq D}\sum_{h\leq H}\frac{1}{h}\sum_{|r|\leq R}\sup_{T\in[N,N+2]}\left|U_{2}^{(2)}\right|.\]

\subsection{\bf The estimation of the sums $\Omega_{1},\Omega_{2}^{(1)},\Omega_{2}^{(2)},\Omega_{3}$}

\begin{lemma}
If $\frac{238}{247}<\gamma<1$, then we have
\[\Omega_{1}\ll\frac{N^{2\gamma-1}}{(\log N)^{2}}, \qquad \Omega_{2}^{(1)}\ll\frac{N^{2\gamma-1}}{(\log N)^{2}}.\]
\end{lemma}
\begin{proof}
In fact, these two formulas are (103) and (102) of \cite{PeTo}, respectively. One can see \cite[Section 3.5]{PeTo} for detail.
\end{proof}

\begin{lemma}
If $\frac{238}{247}<\gamma<1$, then we have
\[\Omega_{2}^{(2)}\ll\frac{N^{2\gamma-1}}{(\log N)^{2}}, \qquad \Omega_{3}\ll\frac{N^{2\gamma-1}}{(\log N)^{2}}.\]
\end{lemma}
\begin{proof}
Let us consider the sum $\Omega_{3}$ first. We divide the sum $U_{3}$ given by \eqref{e.407} into $\mathcal{O}(\log N)$ sums of the form
\begin{equation}\label{e.411}
W_{K,L}=\sum_{L<\ell\leq 2L}b(\ell)\sum_{\substack{K<k\leq 2K \\ \frac{P}{\ell}<k\leq\frac{2P}{\ell}}}a(k)e(h(k\ell)),
\end{equation}
where
\begin{equation}\label{e.412}
a(k)\ll N^{\varepsilon}, \quad b(\ell)\ll N^{\varepsilon}, \quad P^{\frac{1}{3}}\leq K\leq P^{\frac{1}{2}}\ll L\ll P^{\frac{2}{3}}, \quad KL\asymp P,
\end{equation}
and
\begin{equation}\label{e.413}
h(t)=rt^{c}+v\left(T-t^{c}\right)^{\gamma}.
\end{equation}
\vskip 2 pt

It follows from \eqref{e.411}, \eqref{e.412} and Cauchy's inequality that
\begin{equation}\label{e.414}
\left|W_{K,L}\right|^{2}\ll N^{\varepsilon}L\sum_{L<\ell\leq 2L}\left|\sum_{K_{1}<k\leq K_{2}}a(k)e(h(k\ell))\right|^{2},
\end{equation}
where
\[K_{1}=\max\left(K,\frac{P}{\ell}\right), \qquad K_{2}=\min\left(2K,\frac{2P}{\ell}\right).\]
Suppose that $Q$ is an integer which satisfies
\begin{equation}\label{e.415}
1\leq Q\leq K.
\end{equation}
Then by applying Lemma 2.4, we can derive that
\begin{align}
\left|W_{K,L}\right|^{2}&\ll\frac{N^{\varepsilon}LK}{Q}\sum_{L<\ell\leq 2L}\sum_{|q|\leq Q}\left(1-\frac{|q|}{Q}\right) \nonumber\\
& \times\sum_{\substack{K_{1}<k\leq K_{2} \\
K_{1}<k+q\leq K_{2}}}a(k+q)\overline{a(k)}e(h((k+q)\ell)-h(k\ell)) \nonumber\\
& \ll\frac{N^{\varepsilon}LK}{Q}\sum_{|q|\leq Q}\sum_{\substack{K<k\leq 2K \\
K<k+q\leq 2K}}\left|\sum_{L_{1}<\ell\leq L_{2}}e\left(Y_{k,q}(\ell)\right)\right|,\label{e.416}
\end{align}
where
\begin{equation}\label{e.417}
L_{1}=\max\left(L,\frac{P}{k},\frac{P}{k+q}\right), \quad L_{2}=\min\left(2L,\frac{2P}{k},\frac{2P}{k+q}\right)
\end{equation}
and
\begin{equation}\label{e.418}
Y(\ell)=Y_{k,q}(\ell)=h((k+q)\ell)-h(k\ell).
\end{equation}
For $q=0$, by the trivial estimate, we have
\begin{equation}\label{e.419}
\sum_{L_{1}<\ell\leq L_{2}}e\left(Y_{k,q}(\ell)\right)\ll L.
\end{equation}
Hence by \eqref{e.416} and \eqref{e.419}, we get
\begin{equation}\label{e.420}
\left|W_{K,L}\right|^{2}\ll\frac{N^{\varepsilon}(LK)^{2}}{Q}+\frac{N^{\varepsilon} LK}{Q}\sum_{1\leq|q|\leq Q} \sum_{K<k\leq 2K-q}\left|\sum_{L_{1}<\ell\leq L_{2}}e\left(Y_{k, q}(\ell)\right)\right|.
\end{equation}
\vskip 2 pt

Now we consider the function $Y(\ell)$. From \eqref{e.418} we find that
\[Y(\ell)=\int_{k}^{k+q} h_{t}^{\prime}(t \ell) \mathrm{d} t=\int_{k}^{k+q} \ell h^{\prime}(t \ell) \mathrm{d} t\]
and thus
\begin{equation}\label{e.421}
Y^{(j)}(\ell)=\int_{k}^{k+q}\left(jt^{j-1}h^{(j)}(t\ell)+\ell t^{j}h^{(j+1)}(t\ell)\right)\mathrm{d}t, \quad j \geq 1.
\end{equation}
By \eqref{e.413}, \eqref{e.421} and some complicated but elementary calculations, we get
\begin{flalign}
Y'(\ell)&=\int_{k}^{k+q}\left(rc^2(t\ell)^{c-1}+v(t\ell)^{c-1}(T-(t\ell)^c)^{\gamma-2}((t\ell)^c-cT)\right)
\mathrm{d}t\label{e.422}\\
Y^{(j)}(\ell)&=\int_{k}^{k+q}\left(\phi_{1}(t)+\phi_{2}(t)\right)\mathrm{d}t, \qquad j \geq 2,\label{e.423}
\end{flalign}
where
\begin{flalign}
\phi_{1}(t)&=rc^2(c-1)(c-2)\cdots(c-j+1)t^{c-1}\ell^{c-j}, \label{e.424}\\
\phi_{2}(t)&=v(c-1)Tt^{c-1}\ell^{c-j}(T-(t\ell)^c)^{\gamma-j-1}
\sum_{i=1}^{j}g_{i}(c)T^{j-i}(t\ell)^{(i-1)c},\label{e.425}
\end{flalign}
where $g_{i}(x)\in\mathbb{Z}[x]$ are polynomials of degree $j-1$ for $j\geq 2.$
\vskip 2 pt

If $t\in[k,k+q]$ then $t\ell\asymp P$. Thus, by \eqref{e.424}, \eqref{e.425} and the condition $N\leq T\leq N+2$ imposed in \eqref{e.S01}, we find that uniformly for $t\in[k,k+q]$ we have
\begin{equation}\label{e.426}
\left|\phi_{1}(t)\right|\asymp|r|k^{j-1}N^{1-j\gamma} \qquad \text {and} \qquad \phi_{2}(t)\asymp vk^{j-1}N^{-(j-1)\gamma}.
\end{equation}
\vskip 2 pt

From \eqref{e.422}, \eqref{e.423} and \eqref{e.426} we see that there exists a sufficiently small constant $\alpha_{1}>0$ such that if $|r|\leq\alpha_{1}vN^{\gamma-1}$, then $\left|Y^{(j)}(l)\right| \asymp qvk^{j-1}N^{-(j-1)\gamma}$. Similarly, we conclude that there exists a sufficiently large constant $A_{1}>0$ such that if $|r|\geq A_{1}vN^{\gamma-1}$, then $\left|Y^{(j)}(l)\right|\asymp|r|qk^{j-1}N^{1-j\gamma}$. Hence, it makes sense to divide the sum $\Omega_{3}$ into four sums according to the value of $r$ as follows:
\begin{equation}\label{e.427}
\Omega_{3}=\Omega_{3,1}+\Omega_{3,2}+\Omega_{3,3}+\Omega_{3,4},
\end{equation}
where
\begin{flalign}
\text { in } \Omega_{3,1}: & \qquad |r| \leq \alpha_{1} v N^{\gamma-1}, \nonumber\\
\text { in } \Omega_{3,2}: & \qquad -A_{1} v N^{\gamma-1}<r<-\alpha_{1} v N^{\gamma-1}, \nonumber\\
\text { in } \Omega_{3,3}: & \qquad \alpha_{1} v N^{\gamma-1}<r<A_{1} v N^{\gamma-1}, \label{e.4001}\\
\text { in } \Omega_{3,4}: & \qquad A_{1} v N^{\gamma-1} \leq|r| \leq R .\label{e.428}
\end{flalign}
\vskip 2 pt

Let us consider $\Omega_{3,4}$ first. By \eqref{e.409} and \eqref{e.428}, we have
\begin{equation}\label{e.429}
\Omega_{3,4}\ll(\log N)\sum_{d\leq D}\sum_{h\leq H}\frac{1}{h}\sum_{A_{1}vN^{\gamma-1}\leq|r|\leq R}\sup_{T\in[N,N+2]}\left|W_{K,L}\right|.
\end{equation}
\vskip 2 pt
Consider the sum $W_{K,L}.$ We already mentioned that if $|r|\geq A_{1}vN^{\gamma-1}$, then uniformly for $\ell\in(L_{1},L_{2}]$ we have
$$\left|Y^{(j)}(l)\right|\asymp|r|qk^{j-1}N^{1-j\gamma}\asymp |r|qN^{1-\gamma}\cdot(Pk^{-1})^{1-j},\quad j\geq1.$$
Then Lemma 2.3 with $(\kappa,\lambda)=BA^2BA^2(\frac{1}{2},\frac{1}{2})=(\frac{13}{40},\frac{22}{40})$ implies
$$\sum_{L_{1}<\ell\leq L_{2}}e\left(Y_{k,q}(\ell)\right)\ll (|r|qN^{1-\gamma})^{\frac{13}{40}}(Pk^{-1})^{\frac{22}{40}}+(|r|qN^{1-\gamma})^{-1}\ll |r|^{\frac{13}{40}}q^{\frac{13}{40}}k^{-\frac{22}{40}}N^{\frac{13}{40}+\frac{9\gamma}{40}}.$$
Then from \eqref{e.412} and \eqref{e.420} we can see that
\begin{align*}
\left|W_{K,L}\right|^{2}&\ll\frac{N^{\varepsilon}(LK)^{2}}{Q}+\frac{N^{\varepsilon} LK}{Q}\sum_{1\leq|q|\leq Q} \sum_{K<k\leq 2K-q}\left(|r|^{\frac{13}{40}}q^{\frac{13}{40}}k^{-\frac{22}{40}}N^{\frac{13}{40}+\frac{9\gamma}{40}}\right)\\
&\ll N^\varepsilon\left(N^{2\gamma}Q^{-1}+N^\gamma Q^{-1}\sum_{1\leq q\leq Q}\left(|r|^{\frac{13}{40}}q^{\frac{13}{40}}K^{\frac{9}{20}}N^{\frac{13}{40}+\frac{9\gamma}{40}}\right)\right)\\
&\ll N^\varepsilon\left(N^{2\gamma}Q^{-1}+|r|^{\frac{13}{40}}Q^{\frac{13}{40}}N^{\frac{13}{40}+\frac{29\gamma}{20}}\right),
\end{align*}
which implies
\begin{equation}\label{e.430}
W_{K,L}\ll N^\varepsilon\left(N^{\gamma}Q^{-\frac{1}{2}}+|r|^{\frac{13}{80}}Q^{\frac{13}{80}}
N^{\frac{13}{80}+\frac{29\gamma}{40}}\right).
\end{equation}
From \eqref{e.20}, \eqref{e.429} and \eqref{e.430}, we get
\begin{align*}
\Omega_{3,4}&\ll N^{\varepsilon}\sum_{d\leq D}\sum_{h\leq H}\frac{1}{h}\sum_{|r|\leq dN^{1-\gamma}(\log N)^{12}}\left(N^{\gamma} Q^{-\frac{1}{2}}+|r|^{\frac{13}{80}}Q^{\frac{13}{80}}
N^{\frac{13}{80}+\frac{29\gamma}{40}}\right) \\
& \ll N^{\varepsilon}\left(D^{2}NQ^{-\frac{1}{2}}+Q^{\frac{13}{80}}D^{\frac{173}{80}}N^{\frac{53}{40}-\frac{7 \gamma}{16}}\right)\ll N^{\varepsilon}\left(N^{1+2\delta}Q^{-\frac{1}{2}}+Q^{\frac{13}{80}}N^{\frac{53}{40}-\frac{7 \gamma}{16}+\frac{173\delta}{80}}\right).
\end{align*}
Taking
\begin{equation}\label{e.431}
Q=\left[N^{\frac{35\gamma}{53}-\frac{13\delta}{53}-\frac{26}{53}}\right],
\end{equation}
it is easy to see that the condition \eqref{e.415} holds. Hence, for $1<c<\frac{247}{238}$, we get
\begin{equation}\label{e.432}
\Omega_{3,4}\ll\frac{N^{2\gamma-1}}{(\log N)^{2}}.
\end{equation}
\vskip 2 pt

Now let us consider $\Omega_{3,3}$. From \eqref{e.409} and \eqref{e.4001} we have
\begin{equation}\label{e.433}
\Omega_{3,3}\ll(\log N)\sum_{d\leq D}\sum_{h\leq H}\frac{1}{h}\sum_{\alpha_{1}vN^{\gamma-1}<r<A_{1}vN^{\gamma-1}}\sup_{T\in[N,N+2]}\left|W_{K,L}\right|.
\end{equation}
By \cite[(143),(144)]{PeTo}, \eqref{e.417} and \eqref{e.420}, we get
\[W_{K,L}\ll N^{\varepsilon}\left(N^{\gamma}Q^{-\frac{1}{2}}+v^{\frac{1}{4}}Q^{\frac{1}{4}}
N^{\frac{7\gamma}{8}}+v^{-\frac{1}{4}}Q^{-\frac{1}{4}}N^{\frac{7\gamma}{8}}+v^{\frac{1}{12}}Q^{\frac{1}{12}} N^{\frac{11\gamma}{12}}+v^{-\frac{1}{12}}Q^{-\frac{1}{12}}N^{\frac{23\gamma}{24}}\right),\]
Applying the above estimate for $W_{K,L}$ in \eqref{e.433} to get
\begin{align*}
\Omega_{3,3}&\ll N^{\varepsilon}\sum_{d\leq D}\sum_{h\leq H}\frac{1}{h}\sum_{r<A_{1}\log ^{3}N}\left(N^{\gamma} Q^{-\frac{1}{2}}+\left(\frac{h}{d}\right)^{\frac{1}{4}}Q^{\frac{1}{4}}N^{\frac{7\gamma}{8}}\right.\\
&\phantom{=\;\;}\left.+\left(\frac{h}{d}\right)^{-\frac{1}{4}}Q^{-\frac{1}{4}}N^{\frac{7 \gamma}{8}}+\left(\frac{h}{d}\right)^{\frac{1}{12}}Q^{\frac{1}{12}}N^{\frac{11 \gamma}{12}}+\left(\frac{h}{d}\right)^{-\frac{1}{12}}Q^{-\frac{1}{12}}N^{\frac{23\gamma}{24}}\right) \\
&\ll N^{\varepsilon}\left(N^{\gamma+\delta}Q^{-\frac{1}{2}}+Q^{\frac{1}{4}}N^{\frac{1}{4}+\frac{5 \gamma}{8}+\delta}+ Q^{-\frac{1}{4}}N^{\frac{7\gamma}{8}+\frac{5\delta}{4}}+Q^{\frac{1}{12}}N^{\frac{1}{12}+\frac{5 \gamma}{6}+\delta}+Q^{-\frac{1}{12}}N^{\frac{23\gamma}{24}+\frac{13\delta}{12}}\right).
\end{align*}
With the choice of $Q$ which we made in \eqref{e.432} it is now clear that
\begin{equation}\label{e.434}
\Omega_{3,3} \ll \frac{N^{2\gamma-1}}{(\log N)^{2}}.
\end{equation}
\vskip 2 pt
As for $\Omega_{3,1}$ and $\Omega_{3,2}$, by \cite[(146),(147)]{PeTo}, we also have
\begin{equation}\label{e.435}
\Omega_{3,1} \ll \frac{N^{2\gamma-1}}{(\log N)^{2}},\qquad \Omega_{3,2}\ll \frac{N^{2\gamma-1}}{(\log N)^{2}}.
\end{equation}
Then from \eqref{e.427}, \eqref{e.432}, \eqref{e.434} and \eqref{e.435}, we get
$$\Omega_{3}\ll\frac{N^{2\gamma-1}}{(\log N)^{2}}.$$
\vskip 2 pt
To bound $\Omega_{2}^{(2)}$, we can apply the same argument as the one for $\Omega_{3}$ to derive that
$$\Omega_{2}^{(2)}\ll\frac{N^{2\gamma-1}}{(\log N)^{2}}.$$
This completes the proof of Lemma 4.3.
\end{proof}
\vskip 2 pt

Now Proposition 3.1 follows from \eqref{e.408}, \eqref{e.410}, Lemma 4.2 and Lemma 4.3.

\section{\bf Proof of Theorem 1.1}\label{s5}
Now we are in a position to estimate $\Gamma$, which is defined by \eqref{e.301}. Let $\lambda(d)$ be the lower bound Rosser weights of level D. Then by Lemma 2.1, we find
\begin{align}
\Gamma&=\sum_{\substack{P<p\leq 2P, m\in\mathbb{N}\\ \left[p^{c}\right]+\left[m^{c}\right]=N}}(\log p)\sum_{d \mid\left(m,P(z)\right)}\mu(d)\geq\sum_{\substack{P<p\leq 2P, m\in\mathbb{N} \\\left[p^{c}\right]+\left[m^{c}\right]=N}}(\log p)\sum_{d\mid\left(m,P(z)\right)}\lambda(d)\nonumber\\
&=\sum_{d\mid P(z)}\lambda(d)\sum_{P<p\leq 2P}\log p\sum_{\substack{m\in \mathbb{N} \\\left[p^{c}\right]+\left[m^{c}\right]=N\\ m\equiv0(\text{mod}~d)}}1.\label{e.501}
\end{align}
For the innermost sum of \eqref{e.501}, using the trivial identity
$$\sum_{a\leq n<b}1=[-a]-[-b]=b-a-\psi(-b)+\psi(-a)$$
to get
\begin{align}
G_{d,p}&=\sum_{\substack{m\in \mathbb{N} \\\left[p^{c}\right]+\left[m^{c}\right]=N\\ m\equiv0(\text{mod}~d)}}1=
\sum_{\frac{1}{d}\left(N-\left[p^{c}\right]\right)^{\gamma} \leq m<\frac{1}{d}\left(N+1-\left[p^{c}\right]\right)^{\gamma}}1\nonumber\\
&=\frac{\left(N+1-\left[p^{c}\right]\right)^{\gamma}-\left(N-\left[p^{c}\right]\right)^{\gamma}}{d}+
\psi\left(-\frac{1}{d}\left(N-\left[p^{c}\right]\right)^{\gamma}\right)-
\psi\left(-\frac{1}{d}\left(N+1-\left[p^{c}\right]\right)^{\gamma}\right).\label{e.502}
\end{align}
Putting \eqref{e.502} into \eqref{e.501}, we can obtain
\begin{align}
\Gamma&\geq \sum_{d\mid P(z)}\frac{\lambda(d)}{d}\sum_{P<p\leq2P}(\log p)\left(\left(N+1-\left[p^{c}\right]\right)^{\gamma}-\left(N-\left[p^{c}\right]\right)^{\gamma}\right)\nonumber\\
&\quad +\sum_{d|P(z)}\lambda(d)\sum_{P<p\leq2P}(\log p)\left(\psi\left(-\frac{1}{d}\left(N-\left[p^{c}\right]\right)^{\gamma}\right)-
\psi\left(-\frac{1}{d}\left(N+1-\left[p^{c}\right]\right)^{\gamma}\right)\right)\nonumber\\
&=\Gamma_{0}+\Sigma_{0}-\Sigma_{1},\label{e.503}
\end{align}
where
\begin{flalign}
\Gamma_{0}&=\sum_{d\mid P(z)}\frac{\lambda(d)}{d}\sum_{P<p\leq2P}(\log p)\left(\left(N+1-\left[p^{c}\right]\right)^{\gamma}-\left(N-\left[p^{c}\right]\right)^{\gamma}\right)\label{e.504}\\
\Sigma_{j}&=\sum_{d|P(z)}\lambda(d)\sum_{P<p\leq2P}(\log p)\psi\left(-\frac{1}{d}\left(N+j-\left[p^{c}\right]\right)^{\gamma}\right),\quad j=0,1.\label{e.505}
\end{flalign}
\vskip 2 pt
Consider $\Gamma_{0}$. By Chebyshev's prime number theorem and \eqref{e.200}, we get
\begin{align}
A(N)&=\sum_{P<p\leq2P}(\log p)\left(\left(N+1-\left[p^{c}\right]\right)^{\gamma}-\left(N-\left[p^{c}\right]\right)^{\gamma}\right)\nonumber\\
&=\gamma\sum_{P<p\leq2P}(\log p)\left(\left(N-[p^c]\right)^{\gamma-1}+\O(N^{\gamma-2})\right)\asymp N^{2\gamma-1}.\label{e.506}
\end{align}
From \eqref{e.202}, we have
\begin{equation}\label{e.507}
\sum_{d\mid P(z)}\frac{\lambda(d)}{d}\geq \prod_{p<z}\left(1-\frac{1}{p}\right)\left(f(s)+\O((\log D)^{-\frac{1}{3}})\right),
\end{equation}
where
$$s=\frac{\log D}{\log z}=\frac{\delta}{\frac{\delta}{2}-\varepsilon},$$
and where $f(s)$ is given by \eqref{e.02} in Lemma 2.1. Having in mind \eqref{e.02} and \eqref{e.303}, we know that $f(s)>\beta$ for some constant $\beta>0$ depend on $\delta.$ Therefore, by \eqref{e.507} and the Mertens formula we get
$$\sum_{d\mid P(z)}\frac{\lambda(d)}{d}\gg\frac{1}{\log N}.$$
Thus, by \eqref{e.504} and \eqref{e.506} we have
\begin{equation}\label{e.508}
\Gamma_{0}\gg \frac{N^{2\gamma-1}}{\log N}.
\end{equation}
\vskip 2 pt
By combining Proposition 3.1, \eqref{e.503}-\eqref{e.505} and \eqref{e.508}, we finally get
$$\Gamma\gg \frac{N^{2\gamma-1}}{\log N},$$
which proves Theorem 1.1.

\bigskip

\noindent {\bf Acknowledgement.} We wish to thank the referee for a thorough reading of the paper and helpful remarks. The Author would like to express the most sincere gratitude to Professor Yingchun Cai for his valuable advice and
constant encouragement.


\end{document}